\RequirePackage[l2tabu, orthodox]{nag}
\RequirePackage{snapshot}

\documentclass[10pt,twocolumn]{extarticle}
\makeatletter
\if@twocolumn
  \usepackage[dvips,letterpaper,top=0.5in, bottom=0.5in, left=0.75in, right=0.5in,includefoot,heightrounded]{geometry}
\else
  \usepackage[dvips,letterpaper,margin=1in,includefoot,heightrounded]{geometry}
\fi

\usepackage{srcltx}

\usepackage{amsmath}
\usepackage{amssymb,amsfonts}

\usepackage{abstract}

\usepackage{graphicx}
\usepackage{epsfig}
\usepackage[usenames,dvipsnames]{color}
\usepackage[usenames,dvipsnames]{xcolor}
\usepackage{psfrag}

\usepackage{booktabs}

\usepackage{setspace}
\usepackage{flushend}

\usepackage{mwe}

\usepackage{url}

\usepackage{multirow}
\usepackage{algorithmic}
\usepackage[ruled]{algorithm2e}

\usepackage[short,12hr]{datetime}
\usepackage[binary-units]{siunitx}

\usepackage{mathtools}

\DeclarePairedDelimiter{\abs}{\lvert}{\rvert}

\usepackage{hyperref}
\usepackage[hyperpageref]{backref}

\usepackage{fontawesome}
       \usepackage{mathptmx}

\newcommand{\algrule}[1][.2pt]{\par\vskip.5\baselineskip\hrule height #1\par\vskip.5\baselineskip}

\newcommand{\algorithmicinput}{\textbf{Input:}}
\newcommand{\INPUT}{\item[\algorithmicinput]}
\newcommand{\algorithmicoutput}{\textbf{Output:}}
\newcommand{\OUTPUT}{\item[\algorithmicoutput]}

\newcounter{algslikefigs}

\newcommand*\conj[1]{\bar{#1}}

\makeatletter

\newcommand{\printtitle}{%
\makeatletter
\if@twocolumn

\twocolumn[%
  \maketitle
  \begin{onecolabstract}
    \myabstract
  \end{onecolabstract}
  \begin{center}
    \small
    \textbf{Keywords}
    \\\medskip
    \mykeywords
  \end{center}
  \bigskip
]
\saythanks
\else
  \maketitle
  \begin{onecolabstract}
    \myabstract
  \end{onecolabstract}
  \begin{center}
    \small
    \textbf{Keywords}
    \\\medskip
    \mykeywords
  \end{center}
  \bigskip
  \onehalfspacing
\fi
\makeatother
}

\title{%
Fast Matrix Inversion and Determinant Computation for Polarimetric Synthetic Aperture Radar}

\author{%
D. F. G. Coelho%
\thanks{%
D. F. G. Coelho was with the Graduate Program of Electrical and Computer Engineering Department, University of Calgary, Canada.
\mbox{E-mail}: \protect\url{diegofgcoelho@gmail.com}
}
\and
R.~J.~Cintra%
\thanks{%
R. J. Cintra is with the
Signal Processing Group,
Departamento de Estat\'istica,
Universidade Federal de Pernambuco
and
with the
Department of Electrical and Computer Engineering, University of Calgary, Calgary, AB, Canada.
E-mail: \protect\url{rjdsc@de.ufpe.br}
}
\and
A. C. Frery%
\thanks{%
A. C. Frery is with the
Laborat\'orio de Computa\c{c}\~ao Cient\'ifica e An\'alise Num\'erica, Universidade Federal de Alagoas, Brazil.
E-mail: \protect\url{acfrery@gmail.com}
}
\and
V. S. Dimitrov%
\thanks{%
V. S. Dimitrov is with the
Department of Electrical and Computer Engineering, University of Calgary, Calgary, AB, Canada.
}
}

\date{}

\newcommand{\myabstract}{%
This paper
introduces a fast algorithm for simultaneous inversion and determinant computation of small sized matrices in the context of fully Polarimetric Synthetic Aperture Radar (PolSAR) image processing and analysis.
The proposed fast algorithm is based on the computation of the adjoint matrix and the symmetry of the input matrix.
The algorithm is implemented in a general purpose graphical processing unit (GPGPU) and compared to the usual approach based on Cholesky factorization.
The assessment with simulated observations and data from an actual PolSAR sensor show a speedup factor of about~two when compared to the usual Cholesky factorization.
Moreover, the expressions provided here can be implemented in any platform.
}

\newcommand{\mykeywords}{%
PolSAR, Matrix inversion, Determinant, Fast algorithms
}

\begin{document}

\printtitle

\section{Introduction}

Microwave remote sensing is basilar as it provides complementary information to that provided by classical sensors which perceive the optical spectrum.
Longer wavelengths in the
\SI{1}{cm} to the \SI{1}{m}
can penetrate clouds and other adverse atmospheric conditions, as well as canopies and soils.
There are passive and active microwave sensors; the latter carry their own source of illumination, and can be either imaging or non-imaging.
Radar devices belong to the former.
They transmit a radio signal and detect the returned echo (backscatter), with which an image is formed.
Several signal processing techniques are used to enhance the spatial resolution of such imagery.
In particular,
the use of synthetic antennas
leads to synthetic aperture radar (SAR) imaging;
as well as to polarimetric SAR (PolSAR) imaging
which employs polarized signals

The statistical properties of PolSAR were studied in the context of optical polarimetry~\cite{Goodm1985}.
The simplest case occurs when the backscatter is constant.
In such a case, the targets are characterized by a scattering (or Sinclair) matrix~$\mathbf{S}$, which describes dependence of its scattering properties on the polarization.
The scattering matrix is defined on the horizontal~(H) and vertical~(V) basis as
\begin{equation}
\mathbf{S'}=
\left[\begin{array}{cc}
S_{\text{HH}} & S_{\text{HV}}\\
S_{\text{VH}} & S_{\text{VV}}
\end{array}\right],
\label{Eq:scattering_matrix}
\end{equation}
where each element is a complex value, representing the amplitude and the phase of the scattered signal.
Reciprocity~$S_{\text{HV}}=S_{\text{VH}}$ holds for most natural targets, so one may use the scattering
vector~$
\mathbf{S}=
\begin{bmatrix}
S_{\text{HH}} &  S_{\text{HV}} &   S_{\text{VV}}
\end{bmatrix}^\top$
without loss of information,
where~$^\top$ represents transposition~\cite{anfinsen2011goodness}.

More often than not, these single-look complex data are processed in order to improve the signal-to-noise ratio.
Multilook fully polarimetric data are formed as
\begin{equation}
Z^{(N)} = \frac{1}{N} \sum_{\ell=1}^{N} \mathbf{S}(\ell)^{\mathsf{H}}\; \mathbf{S}(\ell),
\label{eq:ZN}
\end{equation}
where ${}^{\mathsf{H}}$ denotes the conjugate transposition and $\ell$ indexes theoretically independent looks of the same scene.
This $3\times3$ complex matrix is positive Hermitian with real entries in the diagonal.

As noted by Torres~\emph{et~al.}~\cite{Torres2014} and the references therein, many techniques for PolSAR image processing and analysis rely on the statistical properties of~$Z^{(N)}$.
The density of several models
(Wishart~\cite{goodm1963statistical},
$K$~\cite{KDistributionPolarimetricPIER90},
$G^0$~\cite{frery1997model},
$G$~\cite{freitas2005polarimetric},
Kummer-U~\cite{doulgeris2015automatic}
to name a few)
depends solely on a few parameters:
the covariance matrix, the number of looks $L$ and, sometimes, texture descriptors.
The covariance matrix is the expected value of $Z^{(N)}$,
namely $\mathbf{Z} = \operatorname{E}\{Z^{(N)}\}$, and it enters the expressions only through its determinant and its inverse.

Moreover, divergences among these models
(Kullback-Leibler, Hellinger, R\'enyi, Bhattacharya, Triangular, Harmonic~\cite{nascimento2010hypothesis}),
test statistics (likelihood ratio)~\cite{conradsen2003test},
and
classification rules~\cite{skriver2012crop}
only require the computation of
$\operatorname{det}(\mathbf{Z})$, $\mathbf{Z}^{-1}$, and~$\operatorname{det}(\mathbf{Z}^{-1})$.
A typical Uninhabited Aerial Vehicle Synthetic Aperture Radar (UAVSAR) image may have $10^{4}\times 4\cdot 10^{4}$ pixels, where each pixel is represented by a~$3\times 3$ Hermitian matrix as in~\eqref{eq:ZN}.
As illustrative example, if one relies on Nonlocal Means Filter approach based on stochastic distances~\cite{Torres2014}, at each sear window, which are typically large, e.g. $23\times 23$ pixels, one needs to compute $23^2$ inverse and $23^2$ determinant operations for each pixel.
This scales to a long computing time as it results in a total of~$23^2\cdot 4\cdot 10^8 \approx 2.1\cdot 10^{11}$ inverse matrix and determinant calculations~\emph{per image}.
Such amount of data is likely to increase with the incoming availability of sensors with finer resolution.
Thus one reaches the conclusion that accurate and fast procedures are of paramount importance for dealing with PolSAR data.

The Cholesky factorization is the most popular numerical analysis method
for the direct solution of linear algebra tasks involving positive definite dense matrices~\cite{Bjoerck2014,Golub1996,Shores2007}.
It is also the algorithm of choice
for matrix inversion and determinant calculation in the context of image classification of PolSAR images~\cite{Torres2014}.
In this paper,
we propose a fast algorithm for the computation of matrix inverse and determinant of~$3\times 3$ Hermitian matrices in the context of PolSAR image classification that outperforms the Cholesky factorization\color{black}.
The introduced algorithm is proven
to reduce the overall arithmetic complexity associated with the matrix inversion and determinant calculation
when compared
with the Cholesky factorization approach.
Such lower arithmetic cost
results in a reduction of the computation time to about a half of the Cholesky based method.
The proposed algorithm and the matrix inversion and determinant calculation based on Cholesky factorization are implemented in  a  general purpose graphical processing unit (GPGPU) using C/C++ and open computing language (OpenCL),
which are tools that have been used for accelerating a several of algorithms
in geoscience and remote sensing~\cite{Lukac2013,Steinbach2012,Li2014,Li2015}.
The proposed algorithm and the method based on Cholesky factorization are tested using
both simulated and
measured PolSAR data~\cite{polsarpro2017}.

The paper unfolds as follows.
Section~\ref{section-mathematical}
introduces notation and preliminary considerations.
The Cholesky factorization is reviewed
and described as a means for
matrix inversion and determinant computation.
In Section~\ref{sec:fast}, the proposed method and its fast algorithm is introduced.
Section~\ref{sec:fast}
brings the arithmetic complexity assessment
and a discussion
on the numerical error analysis
of the proposed method
compared with the Cholesky factorization approach.
Section~\ref{sec:implementation} has implementation considerations including software and hardware comments.
Experiment results shows the effectiveness of the proposed method.
Final comments and suggested directions for future works
are
in
Section~\ref{sec:conclusion}.
\section{Mathematical Review}
\label{section-mathematical}
\subsection{Preliminaries and Notation}
The type of matrix that occurs in the PolSAR problem
is the~$3\times 3$ correlation matrix with complex entries.
Being a correlation matrix,
it inherits the Hermitian property~\cite{Golub1996}.
Therefore,
it can be represented according to:
\begin{align*}
\mathbf{A} =
\begin{bmatrix}
a & b & c\\
\conj{b} & d & e\\
\conj{c} & \conj{e} & f
\end{bmatrix},
\end{align*}
where
$a$, $d$, and $f$
are real quantities;
$b$, $c$, and~$e$ are complex numbers;
and
the overbar $\ \conj{}\ $ denotes the complex conjugation.
Because~$\mathbf{A}$ is Hermitian,
in actual implementations,
only the upper or lower triangular part
of the matrix is stored for computation.

The~$(i,j)$ cofactor of the matrix~$\mathbf{A}$ is the determinant of the submatrix formed with the elimination of the~$i$th row and~$j$th column times~$(-1)^{i+j}$~\cite{Golub1996}.
The adjoint matrix is computed according to the following~\cite{Golub1996}:
\begin{align*}
\tilde{\mathbf{A}} =
\begin{bmatrix}
a_c & b_c & c_c\\
\conj{b}_c & d_c & e_c\\
\conj{c}_c & \conj{e}_c & f_c
\end{bmatrix},
\end{align*}
where the element in position~$(i,j)$ represents the cofactor of the element~$(j,i)$ in the original matrix~$\mathbf{A}$  and are given by
\begin{equation}
\label{eq:adjoint_elements}
\begin{split}
a_c & = d\cdot f-\abs{e}^2
,
\\
b_c & = c\cdot \conj{e}-b\cdot f
,
\\
c_c & = b\cdot e-c\cdot d
,
\\
d_c & = a\cdot f-\abs{c}^2
,
\\
e_c & = c\cdot \conj{b}-a\cdot e
,
\\
f_c & = a\cdot d-\abs{b}^2
.
\end{split}
\end{equation}
The adjoint matrix~$\tilde{\mathbf{A}}$ is also Hermitian;
thus~$a_c, d_c,$ and~$f_c$ are real numbers.

\subsection{Cholesky Factorization}
\label{subsec:cholesky}

Cholesky factorization is applied in many numerical problems~\cite{Aquilante2008,Wilson1990,Kershaw1978}.
Let~$\mathbf{A}$ be the input matrix we are interested into inverting and computing the determinant.
Cholesky factorization decomposes an input matrix
into the product
$\mathbf{L} \cdot \mathbf{L}^{\mathsf{H}}$,
where
$\mathbf{L}$ is a lower triangular matrix
and ${}^\mathsf{H}$ represents the transpose conjugate operation.
Let~$l_{i,j}$ be the~$(i,j)$ entry of~$\mathbf{L}$.
The Cholesky factorization is based on the
following relations
between the elements of~$\mathbf{A}$ and~$\mathbf{L}$:
\begin{equation*}
\begin{split}
\abs{l_{0,0}}^2 & = a
,
\\
l_{1,0} \cdot l_{0,0} & = \conj{b}
,
\\
l_{2,0} \cdot l_{0,0} & = \conj{c}
,
\\
\abs{l_{1,1}}^2+\abs{l_{1,0}}^2 & = d
,
\\
l_{2,1} \cdot l_{1,1}+l_{1,0} \cdot \conj{l}_{2,0} & = e
,
\\
\abs{l_{2,2}}^2+\abs{l_{2,0}}^2+\abs{l_{2,1}}^2 & = f
,
\end{split}
\end{equation*}
where~$\abs{\cdot}$
returns the magnitude of its complex argument~\cite{Graham1989}.

Once matrix $\mathbf{L}$ is derived,
its inverse can be directly
obtained from the following expressions:
\begin{align*}
\mathbf{L}^{-1} =
\begin{bmatrix}
\frac{1}{l_{0,0}} & 0 & 0\\
\frac{1}{l_{1,1}}\cdot \left(1-\frac{l_{1,0}}{l_{0,0}}\right) & \frac{1}{l_{1,1}} & 0\\
\frac{1}{l_{2,2}}
\cdot
\left(
\frac{l_{2,1}\cdot l_{1,0}}{l_{1,1}\cdot l_{0,0}}
-\frac{l_{2,0}}{l_{0,0}}
\right)
&
-\frac{l_{2,1}}{l_{1,1}}\frac{1}{l_{2,2}} & \frac{1}{l_{2,2}}
\end{bmatrix}
\end{align*}
Given~$\mathbf{L}^{-1}$,
the inverse of the input matrix~$\mathbf{A}$
is furnished by:
$\mathbf{A}^{-1}
=
{(\mathbf{L}^{-1})}^\mathsf{H}
\cdot
\mathbf{L}^{-1}$.
The determinant
of $\mathbf{A}$,
$\operatorname{det}(\mathbf{A})$,
is also available as a by-product
of the discussed factorization~\cite{Golub1996}.
The detailed procedure
for inverting~$3\times 3$ matrices
based on Cholesky factorization
is described in~Algorithm~\ref{alg:cholesky}.

\begin{figure}
\hrule
\begin{algorithmic}
\medskip

\INPUT
$a, d, f \in \mathbb{R}$ and~$b, c, e \in \mathbb{C}$

\OUTPUT
$a_i, d_i, f_i \in \mathbb{R}$ and~$b_i, c_i, e_i \in \mathbb{C}$ and~$\operatorname{det}(\mathbf{A})$

\algrule

\STATE Compute~$\mathbf{L}$ and auxiliary variables:\\
\begin{tabular}{l l}
$l_{0,0} \gets \sqrt{a}$ & $v_0 \gets 1/l_{0,0}$\\
$l_{1,0} \gets \conj{b}\cdot v_0$ & $l_{2,0} \gets \conj{c}\cdot v_0$\\
$l_{1,1} \gets \sqrt{d-\abs{l_{1,0}}^2}$ & $v_1\gets 1/l_{1,1}$\\
$l_{2,1} \gets (e-l_{1,0}\cdot \conj{l}_{2,0})\cdot v_1$ &
$l_{2,2} \gets \sqrt{f-\abs{l_{2,0}}^2-\abs{l_{2,1}}^2}$\\
$v_2 \gets 1/l_{2,2}$
\end{tabular}

\algrule

\STATE Compute~$\mathbf{L}^{-1} = \{ t_{i,j} \}$:\\

\begin{tabular}{l l}
$t_{0,0} \gets
v_0$
&
$t_{1,0} \gets
-l_{1,0}\cdot v_0\cdot v_1$
\\
$t_{1,1} \gets
v_1$
&
$t_{2,0} \gets
v_2\cdot v_0\cdot \left(l_{2,1}\cdot l_{1,0}\cdot v_1-l_{2,0}\right)$
\\
$t_{2,1} \gets
-l_{2,1}\cdot v_1\cdot v_2$
&
$t_{2,2} \gets
v_2$
\end{tabular}

\algrule

\STATE Compute the determinant:
$\operatorname{det}(\mathbf{A}) \gets
(l_{0,0}\cdot l_{1,1}\cdot l_{2,2})^2$

\algrule

\STATE Compute~$\mathbf{A}^{-1}$:\\
\begin{tabular}{l l}
$a_i \gets
t_{0,0}^2
+
\abs{t_{1,0}}^2
+
\abs{t_{2,0}}^2$
&
$b_i \gets
\conj{t}_{1,0}\cdot t_{1,1}
+
\conj{t}_{2,0}\cdot t_{2,1}$
\\
$c_i \gets
\conj{t}_{2,0} \cdot t_{2,2}$
&
$d_i \gets
t_{1,1}^2+\abs{t_{2,1}}^2$
\\
$e_i \gets
\conj{t}_{2,1} \cdot t_{2,2}$
&
$f_i \gets
t_{2,2}^2$
\end{tabular}

\algrule

\RETURN{$a_i, b_i, c_i, d_i, e_i, f_i$ and~$\operatorname{det}(\mathbf{A})$}

\end{algorithmic}
\hrule

\stepcounter{algslikefigs}
\renewcommand{\thefigure}{\arabic{algslikefigs}}
\renewcommand{\figurename}{Alg.}

\caption{The matrix inversion and determinant calculation based on
the Cholesky factorization.}
\label{alg:cholesky}

\end{figure}

\section{Fast Inversion and Determinant Computation}
\label{sec:fast}

\subsection{Proposed Fast Algorithm}

From matrix theory~\cite[p.~59]{seber2008matrix},
we know that
\begin{equation}
\label{eq:inv_adjoint}
\mathbf{A}^{-1}
=
\frac{1}{\operatorname{det}(\mathbf{A})}\tilde{\mathbf{A}}
\end{equation}
and
\begin{equation}
\label{eq:det}
\operatorname{det}(\mathbf{A})
=
a\cdot a_c
+
\conj{b}\cdot b_c
+
\conj{c}\cdot c_c
.
\end{equation}
Notice that~\eqref{eq:det} is a direct consequence of the fact that the determinant of a matrix can be obtained by the summation of any of the row or column elements weighted by their corresponding cofactors~\cite{Watkins2004}.
Notice that the first term in~\eqref{eq:det},
$aa_c$,
is purely real.
Also notice that because $\mathbf{A}$ is Hermitian,
its determinant is real-valued~\cite{Golub1996}.
Therefore,
the term~$(\conj{b}\cdot b_c+\conj{c}\cdot c_c)$ must be real.
As a consequence,
we calculate
only the real part of~$\conj{b}b_c+\conj{c}c_c$ in order to compute~$\operatorname{det}(\mathbf{A})$
as follows:
\begin{equation}
\label{eq:detr}
\begin{split}
\operatorname{det}(\mathbf{A})
=
&
a\cdot a_c
+
\Re{(b)}\cdot \Re{(b_c)}
+
\Im{(b)}\cdot \Im{(b_c)}
\\
&
+
\Re{(c)}\cdot \Re{(c_c)}
+
\Im{(c)}\cdot \Im{(c_c)}
,
\end{split}
\end{equation}
where~$\Re{(\cdot)}$ and~$\Im{(\cdot)}$ return
the real and imaginary parts of its complex argument,
respectively.

Expressions~\eqref{eq:adjoint_elements},
\eqref{eq:inv_adjoint},
and~\eqref{eq:detr}
can be combined
to derive a fast algorithm
for the computation of~$\mathbf{A}^{-1}$
and~$\operatorname{det}(\mathbf{A})$.
Algorithm~\ref{alg:fast}
details the efficient procedure.
As input data,
it requires only the upper triangular part of the input matrix~$\mathbf{A}$.
The proposed algorithm
returns
(i)~the upper triangular part of $\mathbf{A}^{-1}$,
which can be fully reconstructed due to the Hermitian property
and
(ii)~inverse and determinant.
The quantity~$1/\operatorname{det}(\mathbf{A})$,
which appears in all the aforementioned densities and divergences,
is readily available
as a by-product of Algorithm~\ref{alg:fast}
under the auxiliary variable~$t$.
Notice that, although not explicitly stated,
both Algorithm~\ref{alg:cholesky} and Algorithm~\ref{alg:fast}
can be adapted
to check
whether
the matrices subject to calculation are full-rank.
This can be done by simply checking if~$\operatorname{det}(\mathbf{A}) \neq 0$.
This does not represent any significant overhead since~$\operatorname{det}(\mathbf{A})$ is actually necessary for the computation.

\begin{figure}

\hrule
\begin{algorithmic}
\medskip

\INPUT
$a, d, f \in \mathbb{R}$ and~$b, c, e \in \mathbb{C}$

\OUTPUT
$a_i, d_i, f_i \in \mathbb{R}$ and~$b_i, c_i, e_i \in \mathbb{C}$ and~$\operatorname{det}(\mathbf{A})$

\algrule

\STATE Compute cofactors:\\
\begin{tabular}{l l}
$a_i \gets d\cdot f-\abs{e}^2$ &
$b_i \gets c\cdot \conj{e}-b\cdot f$\\
$c_i \gets b\cdot e-c\cdot d$ &
$d_i \gets a\cdot f-\abs{c}^2$\\
$e_i \gets c\cdot \conj{b}-a\cdot e$ &
$f_i \gets a\cdot d-\abs{b}^2$
\end{tabular}
\algrule
\STATE Compute determinant:
\begin{flalign*}
\operatorname{det}(\mathbf{A}) &
\gets
a\cdot a_i
+
\Re{(b)}\cdot \Re{(b_i)}\phantom{=}+
\Im{(b)}\cdot \Im{(b_i)}
\\&\phantom{=}+
\Re{(c)}\cdot \Re{(c_i)}
+
\Im{(c)}\cdot \Im{(c_i)}&&
\end{flalign*}
\algrule
\STATE Compute the inverse of the determinant:
$t \gets 1/\operatorname{det}(\mathbf{A})$

\algrule
\STATE Compute the~$\mathbf{A}^{-1}$ entries:\\
\begin{tabular}{l l}
$a_i \gets t\cdot a_i$ &  $b_i \gets t\cdot b_i$\\
$c_i \gets t\cdot c_i$ &  $d_i \gets t\cdot d_i$\\
$e_i \gets t\cdot e_i$ &  $f_i \gets t\cdot f_i$
\end{tabular}

\algrule

\RETURN{$a_i, b_i, c_i, d_i, e_i, f_i$ and~$\operatorname{det}(\mathbf{A})$}

\end{algorithmic}
\hrule

\stepcounter{algslikefigs}
\renewcommand{\thefigure}{\arabic{algslikefigs}}
\renewcommand{\figurename}{Alg.}

\caption{Proposed algorithm for fast matrix inversion and determinant calculation.}
\label{alg:fast}

\end{figure}

\subsection{Complexity Assessment and Discussion}
\label{sec:complexity}

Algorithm~\ref{alg:cholesky} requires the computation of the squared norm of six complex quantities.
This demands~12 multiplications
and six additions of real quantities~\cite{Blahut2010}.
Besides these operations,
Algorithm~\ref{alg:cholesky} demands~44 multiplications,
12~additions,
three inversions,
three squarings,
and three square roots of real quantities.
This accounts for a total of~89 operations
as shown in Table~\ref{tab:complexity}.
Note that the proposed Algorithm~\ref{alg:fast} outperforms
Algorithm~\ref{alg:cholesky} not only in the total count of operations.
In fact, it reduces the arithmetic cost for all
operations listed in Table~\ref{tab:complexity},
except for the number of additions which requires only two extra additions.

\begin{table}
\centering
\caption{Operation Counts
in Terms of Real Arithmetic Operations with Real Inputs~$x$ and~$y$}
\begin{tabular}{c@{\quad}c@{\quad}c@{\quad}c@{\quad}c@{\quad}c@{\quad}}\toprule
Method & $x\cdot y$ & $x+y$ & $1/x$ & $\sqrt{x}$ & Total
\\
\midrule
Cholesky &          59 & 24 & 3 & 3 & 89 \\
\textbf{Proposed} & 37 & 26 & 1 & 0 & 64 \\
\bottomrule
\end{tabular}
\label{tab:complexity}
\end{table}

In contrast,
the proposed
Algorithm~\ref{alg:fast} requires~four multiplications and two additions for the computation of~$a_i$, $d_i,$ and~$f_i$, each.
It also requires ten multiplications and six additions for the computation of~$b_i$, $c_i,$ and~$e_i$, each.
The computation of~$\operatorname{det}(\mathbf{A})$ and its inverse demands five multiplications, four additions, and only one inversion.
Because of that, Algorithm~\ref{alg:fast} requires a total of~37 multiplications, 26~additions, and one division of real quantities.
A total of 64~operations; 28\% less than
the usual Cholesky factorization approach.
Above quantities are summarized in Table~\ref{tab:complexity}.

\subsection{Numerical Errors}
\label{sec:errors}

Modern arithmetic logic units (ALUs) in contemporary processors
are built to provide high accuracy to a wide range of arithmetic functions.
Under certain conditions,
the accuracy obtained
for
fundamental operations,
such as addition and multiplication,
can also be achieved
for other arithmetic functions,
such as division and square-root operations.
This is not obvious since more complex operations
require a number of additions and multiplications
and are usually implemented through iterative methods~\cite{Kulisch1981}.
However,
in modern GPGPUs
that
follow the OpenCL and CUDA standards,
such as the one employed in our computations,
division and square-root operations
are correctly rounded up to the very
last bit
as it is also the case for
addition and multiplication~\cite[cf. Table 7]{CUDA2018}~\cite[cf. Table 7.2]{Munshi2012}.

The only operations required by the Algorithm~\ref{alg:cholesky} and Algorithm~\ref{alg:fast} are addition, multiplication, division, and square-root operations.
Because all those operations are computed with the same precision,
their numerical accuracies
are the result of the amount of operations required by each of them.
In particular,
quantities that require the largest number of operations
are expected to
dominate the
the numerical error.
This is given by the critical path in terms of the number of arithmetic operations.
A careful analysis of Algorithm~\ref{alg:cholesky}
shows that one needs 37~operations along the critical path to obtain the desired output.
There are more than one critical path,
but as an example,
one of them is:~$b_i \gets t_{2,0}\gets l_{2,1}\gets l_{2,0}\gets v_0 \gets l_{0,0}\gets a$.
In the case of Algorithm~\ref{alg:fast},
the critical path requires 22~operations.
A possible path is:~$b_i \gets t \gets \operatorname{det}(\mathbf{A}) \gets b_i \gets b$.
Because of that, the proposed method in Algorithm~\ref{alg:fast}
can compute the matrix inverse and determinant
with
a higher accuracy
when compared with
the method in Algorithm~\ref{alg:cholesky}.

\section{Implementation and Results}
\label{sec:implementation}

\subsection{Material and Equipment}

We implemented the proposed Algorithm~\ref{alg:fast} and the Cholesky factorization based approach
using open computing language (OpenCL) with
the C~application interface (API)~\cite{AMD2010};
being
compiled for a general purpose graphical processing unit (GPGPU).
The device used is a NVIDIA GeForce 940M, which contains~three computing units and allows the data to be scheduled through the maximum of three different dimensions.
Each dimension allows the maximum use of 1024, 1024, and 64 work items at the same time, respectively.
The device has a global memory of~\SI{2}{\giga\byte}.
The host side was implemented in C/C++ language.
The machine used has a Core~i7 microprocessor and \SI{8}{\giga\byte} of random access memory (RAM) and runs a Linux operating system with distribution Ubuntu version~16.04~LTS.

The communication between the proposed method
and popular platforms for
image processing, such as Matlab, ENVI, and NEST
can be performed by means of interface capabilities offered by each different vendor.
For Matlab, this can be accomplished by mex functions.
ENVI or NEST
can resort to
direct calls of the
compiled C/C++ code
to
interface with
the proposed algorithm.

\subsection{Data Setup}

We separated simulated and measured data
for numerical analysis.

\subsubsection{Simulated Data}

Simulated data were generated from random matrices
obtained with the Eigen C/C++ library~\cite{Guennebaud2010}
 through the method~\texttt{Eigen::MatrixXcd::Random}.
The real and imaginary parts of the matrix coefficients are drawn from the uniform distribution on the closed interval~$[-1,1]$.
The obtained matrices are then converted into Hermitian matrices by computing the product of itself with the transpose conjugate~\cite{Golub1996}.

\subsubsection{Measured Data}

For the measured data,
we considered the data
available in~\cite{polsarpro2017}.
Such data were acquired by the Airborne Synthetic Aperture Radar (AIRSAR) sensor over San Francisco, CA, in
fully polarimetric mode and L-band.
The approximate spatial resolution
is $\SI{10}{m}\times\SI{10}{m}$.
The data were stored in~\texttt{.rdata}
format, original from R language~\cite{RDCT2008}.
For recovering the data, we employ RInside C/C++ library~\cite{Eddelbuettel2013}, which allows the binding of R and C/C++ applications.
The data were converted to the appropriate types in C/C++ and then sent to the GPGPU using OpenCL scheduling runtime routines.
The data used for this experiment has a total of~$450\times 600$ pixels.

\subsection{Methodology and Results}

The
OpenCL
implementations of Algorithm~\ref{alg:cholesky} and Algorithm~\ref{alg:fast}
were executed using all the compute
units and their respective work items
as a single dimension~\cite[pg. 41]{AMD2010}\cite{AMD2015}.
The GPGPU device was allowed to manage and use the preferred
work-group size~\cite{AMD2010}.
The computation
was performed
with \num{1000}~replicates in order to reduce the operating system
fluctuations when estimating the speedup.

For the simulated data computation,
we considered
\num{5469354}~matrices for each replicate.
This amount of matrices
was carefully selected
to
use the GPGPU device
at the limit imposed by
its available memory for allocation,
and corresponds to an image of approximately $2340\times2340$~pixels.
Indeed,
GPGPU finite registers (buffers)
are dedicated to data interchange with the host device
and
the buffer size sets an upper bound to the maximum allocated memory.
This is a limiting factor in the high performance applications.
For measured data,
we submitted
the \num{270000}~matrices
associated
to the selected 450$\times$600 pixel image.

Table~\ref{table-results}
displays
descriptive statistics
of the computing time
required by Algorithm~\ref{alg:cholesky} and Algorithm~\ref{alg:fast}
for
simulated
and
measured
data.
We also provide the ratio between
the measures
as figure of merit
for the speed gain.
Figure~\ref{fig:boxplot}
displays
boxplots of the runtimes for the Cholesky based approach and the proposed method both for measured and simulated data.
The time statistics for simulated data are
larger when compared to the results from measured data,
because
simulated and measured data sets
have
different number of matrices.
However,
the ratios of the computation time statistics
for the simulated and measure data
are very close.

\begin{table*}

\centering

\caption{%
Run-time execution statistics:
minimum ($t_\text{min}$),
average ($t_\text{avg}$),
maximum ($t_\text{max}$),
and standard deviation ($t_\text{std}$)
for matrix inversion and determinant evaluation.
Measurements in milliseconds}

\label{table-results}

\begin{tabular}
{c@{\quad}c@{\quad}c@{\quad}c@{\quad}c@{\quad}c@{\quad}c@{\quad}}
\toprule
\multirow{2}{*}{Time Statistics} &
\multicolumn{3}{c}{Simulated Data (\num[scientific-notation = true,round-mode = places,round-precision = 1]{5469354} matrices)} &
\multicolumn{3}{c}{Measured Data (\num[scientific-notation = true,round-mode = places,round-precision = 1]{270000} matrices)}
\\
\cmidrule(r){2-4}
\cmidrule(r){5-7}
&
Cholesky &
\textbf{Proposed} &
Gain &
Cholesky &
\textbf{Proposed} &
Gain
\\
\midrule
$t_\text{min}$ & 498.89 & 223.06 & 2.23 & 24.68 & 11.18 & 2.21 \\
$t_\text{avg}$ & 507.54 & 231.51 & 2.19 & 28.35 & 12.98 & 2.18 \\
$t_\text{max}$ & 556.79 & 273.08 & 2.03 & 62.60 & 28.59 & 2.18 \\
$t_\text{std}$ &   9.34 &   6.11 & --   &  2.55 &  1.44 & --   \\
\bottomrule
\end{tabular}

\end{table*}

\setcounter{figure}{0}
\begin{figure}
\centering
\psfrag{Cholesky}[][][0.75]{Cholesky}
\psfrag{Proposed}[][][0.75]{Proposed}
\psfrag{40}[][][0.75]{40}
\psfrag{30}[][][0.75]{30}
\psfrag{25}[][][0.75]{25}
\psfrag{20}[][][0.75]{20}
\psfrag{15}[][][0.75]{15}
\psfrag{10}[][][0.75]{10}
\psfrag{5}[][][0.75]{5}
\psfrag{0}[][][0.75]{0}

\psfrag{150}[][][0.75]{150}
\psfrag{300}[][][0.75]{300}
\psfrag{450}[][][0.75]{450}
\psfrag{600}[][][0.75]{600}

\psfrag{T}[][][0.75]{$t$ (ms)}

\epsfig{file=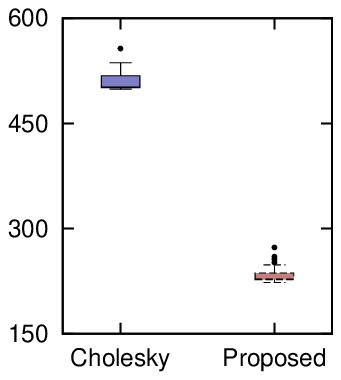, width=4.3cm, height=3.5cm}
\epsfig{file=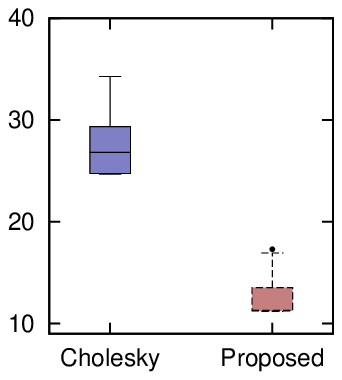, width=4.3cm, height=3.5cm}

\caption{Boxplots for the computing time (in milliseconds)
of the Cholesky based approach and the proposed method using
simulated data (left)
and
measured data (right).
}
\label{fig:boxplot}

\end{figure}

One may notice that PolSAR images can be larger than the images considered adopted in this work.
Indeed,
images to be processed may not fit
into the global memory
of a particular GPGPU adopted for the data processing.
In this case,
the data can be streamed into the GPU
to achieve its maximum capability.
One can start from the left-top corner of the image and stream the data as it is being fed into the GPGPU.
After the maximum amount of data is sent, the data can be processed and then pulled out.
The process must be repeated
until all images have been processed.
The maximum capability, naturally,
depends on the device specifications
of
the computational platform
employed for processing the PolSAR images.

\subsection{Discussion}

Although Table~\ref{tab:complexity}
suggests an improvement of~$89/64 \approx 1.4$,
i.e., roughly a $28\%$ reduction in complexity,
the actual improvement gain or complexity savings is much larger.
Indeed,
to correctly evaluate the speed improvements,
one must take into account that inversion and square-root
computations
demand more GPU clock cycles
when compared with
simple additions or multiplications.
In particular, reducing the number of calls for
the square root operations proves to be significant.
This is because,
in modern processors,
the square root operation
is not computed directly
but through the calculation of the inverse square root~\cite{Parhami2010},
which requires calling a division operation
at the end of the calculation
in order to obtain
the sought result~\cite{Hasnat2017, Libessart2017, Karp1997}.
Therefore,
the total impact of the arithmetic reduction in the proposed algorithm is not just the ratio of the counting of arithmetic operations of the Cholesky based approach and the proposed fast algorithm.
Rather, this ratio is a conservative lower bound.
Results from Table~\ref{table-results}
confirm this analysis,
showing that the speedup gain is more than double.

\section{Final Comments and Future Works}
\label{sec:conclusion}

\color{black}

In this paper, we reviewed the Cholesky factorization method for the inversion and determinant calculation of~$3 \times 3$ Hermitian matrices in the context of PolSAR image processing.
As the main contribution, we proposed a fast algorithm for the computation of matrix inverse and determinant computation for such matrices.
The proposed algorithm is based on the computation of the adjoint matrix and careful examination of the quantities involved.
The derived algorithm was compared to the usual Cholesky based factorization approach and shown to significantly reduce the total number of arithmetic operations.
The proposed algorithm and the Cholesky based approach were both implemented in GPGPU using OpenCL and C/C++.
The proposed approach proved to reduce the total computation time in about a half when compared to the Cholesky based approach with real data;
in other words, the proposed method offers
a 120.7\% improvement in speed, at least.

Future works may include the combination of several
matrix inversions into
a single matrix inversion call
and the application to multifrequency PolSAR images~\cite{Frery2007}.
In this particular case, the matrices to be inverted
and have their
determinants calculated follow
a block diagonal structure,
where each block is
also Hermitian.
If the images are obtained with a single look~$L = 1$,
the inverse and the determinant calculation of the
matrices of interest can be
obtained after applications
of the proposed
method to each of the diagonal blocks.

\section*{Acknowledgments}
Authors acknowledge
CNPq,
FAPEAL,
and NSERC
for the partial support.

{\small
\singlespacing
\bibliographystyle{ieeetr}
\bibliography{bibcleanoutput}
}

\end{document}